\documentclass[MSNbibl,number,citesort,seceqn,dvips]{arxbj}
\usepackage{upgreek}
\usepackage{graphicx}

\aid{0}
\volume{18}
\issue{4}
\pubyear{2012}
\firstpage{1249}
\lastpage{1266}
\doi{10.3150/12-BEJ378} 

\makeatletter
\newtheorem{them}{Theorem}[section]
\newremark{rem}{Remark}[section]
\newtheorem{proposition}{Proposition}[section]
\newtheorem{corollary}{Corollary}[section]

\makeatother

\begin{document}
\begin{frontmatter}

\title{Bootstrap confidence intervals for isotonic estimators in a stereological problem}
\runtitle{Bootstrap in a stereological problem}

\begin{aug}
\author[1]{\fnms{Bodhisattva} \snm{Sen}\corref{}\thanksref{1}\ead[label=e1]{bodhi@stat.columbia.edu}} \and
\author[2]{\fnms{Michael} \snm{Woodroofe}\thanksref{2}\ead[label=e2]{michaelw@umich.edu}}
\runauthor{B. Sen and M. Woodroofe} 
\address[1]{Department of Statistics, Columbia University, New York, NY
10027, USA.\\ \printead{e1}}
\address[2]{Department of Statistics, University of Michigan, Ann
Arbor, MI 48109, USA.\\ \printead{e2}}
\end{aug}

\received{\smonth{9} \syear{2010}}
\revised{\smonth{4} \syear{2011}}

%
\begin{abstract}
Let $\mathbf{X} = (X_1,X_2,X_3)$ be a spherically symmetric random
vector of which only $(X_1,X_2)$ can be observed. We focus attention
on estimating $F$, the distribution function of the squared radius
$Z := X_1^2 + X_2^2 + X_3^2$, from a random sample of $(X_1,X_2)$.
Such a problem arises in astronomy where $(X_1,X_2,X_3)$ denotes the
three dimensional position of a star in a galaxy but we can only observe
the projected stellar positions $(X_1,X_2)$. We consider isotonic
estimators of $F$ and derive their limit distributions. The results are
nonstandard with a rate of convergence $\sqrt{n/{\log n}}$. The
isotonized estimators of $F$ have exactly half the limiting variance
when compared to naive estimators, which do not incorporate the shape
constraint. We consider the problem of constructing point-wise confidence
intervals for $F$, state sufficient conditions for the consistency of a
bootstrap procedure, and show that the conditions are met by the
conventional bootstrap method (generating samples from the empirical
distribution function).
\end{abstract}

%

\end{frontmatter}

\section{Introduction}
Stereology is the study of three-dimensional properties of objects or
matter usually observed two-dimensionally. We consider such a problem,
which arises in astronomy. Suppose that the position $\mathbf{X} :=
(X_1,X_2,X_3)$ of a star within a given galaxy has a spherically
symmetric distribution and that we observe the projected stellar
positions, that is, $(X_1,X_2)$ (with a proper choice of co-ordinates);
and consider the problem of estimating the distribution function $F$ of
the squared distance $Z := X_1^2 + X_2^2 + X_3^2$ of a star to the
center of the galaxy from a random sample of $(X_1,X_2)$. In this
paper, we study the statistical properties of three estimators of $F$.
We show that enforcing \textit{known} shape restrictions (monotonicity)
in the estimation procedure leads to estimators with lower asymptotic
variance (exactly by one-half in this case). We also consider the
problem of constructing point-wise confidence intervals (CIs) around
$F$, and show that the conventional bootstrap method can be used to
construct valid CIs. Our treatment is similar in flavor to Groeneboom and Jongbloed's
\cite{gj95} study of the Wicksell's \cite{w25} ``Corpuscle Problem.''

Suppose that $\mathbf{X}$ has a density of the form $\rho(x_1^2 +x_2^2
+ x_3^2)$. Then $Y := X_1^2 + X_2^2 \sim G$ and $Z$ have densities
\begin{eqnarray}\label{eq:g}
g(y) = \uppi\int_{y}^{\infty} \frac{\rho(z)}{\sqrt{z - y}} \,\mathrm{d}z
\end{eqnarray}
and $f(z) = 2\uppi\sqrt{z}\rho(z)$. The reader may recognize (\ref{eq:g})
as Abel's transformation. It may be inverted as follows. Let
\begin{eqnarray*}
V(y) = \int_y^{\infty}\frac{g(u)}{\sqrt{u - y}} \,\mathrm{d}u.
\end{eqnarray*}
Then
\begin{eqnarray}\label{eq:SimplifyV}
V(y) = \uppi\int_y^{\infty} \biggl[\int_{u}^{\infty} {\rho(z)\,\mathrm{d}z\over\sqrt
{z-u}} \biggr]{\,\mathrm{d}u\over\sqrt{u-y}} = \uppi^2\int_{y}^{\infty} \rho(z)\,\mathrm{d}z
\end{eqnarray}
so that $\rho(z) = -V'(z)/\uppi^2$ at continuity points. Observe that $V$
is a \textit{nonincreasing} function. The quantity of interest, $F$,
can be related to $V$ and, therefore, to the distribution of
$(X_1,X_2)$ by
\begin{eqnarray}\label{eq:F}
F(x) = \int_0^x 2 \uppi\sqrt{u} \rho(u) \,\mathrm{d}u = 1 + \frac{2}{\uppi} \int
_x^\infty\sqrt{z} \,\mathrm{d}V(z),
\end{eqnarray}
where the last equality follows from $\int_0^\infty 2 \uppi\sqrt{u}\rho
(u) \,\mathrm{d}u = 1$. Relationship (\ref{eq:F}) will be used extensively in the
sequel. Let
\begin{eqnarray}\label{eq:V_U}
U(x) := \int_{0}^{x} V(t) \,\mathrm{d}t \nonumber
\end{eqnarray}
for $x > 0$. Then $U$ is concave since $V$ is nonincreasing. Concavity
can also be seen from
\begin{eqnarray}\label{eq:Simplify_U}
U(x)  =  2 \int_0^{\infty} \bigl\{\sqrt{u} - \sqrt{(u - x)}_+\bigr\} g(u) \,\mathrm{d}u,
\nonumber
\end{eqnarray}
where $y_+ = \max\{y,0\}$. Let $J(t) := \int_t^\infty\sqrt{z - t}
\,\mathrm{d}V(z)$. Then
\begin{eqnarray}\label{eq:Simplify_G}
G(t) &=& \uppi\int_{0}^{\infty} \int_{0}^{t \wedge z} \frac{\rho
(z)}{\sqrt{z - y}} \,\mathrm{d}y \,\mathrm{d}z \nonumber
\\[-8pt]
\\[-8pt]
& = & 2 \uppi\int_0^{\infty} \bigl\{\sqrt{z} - \sqrt{(z - t)_+}\bigr\} \rho(z) \,\mathrm{d}z
= 1 + \frac{2}{\uppi} J(t),\nonumber
\end{eqnarray}
where the last step follows from $\int_{0}^{\infty} 2\uppi\sqrt{z}\rho
(z)\,\mathrm{d}z = 1$ and $J(t) = - \uppi^2 \int_t^\infty\sqrt{z - t} \rho(z) \,\mathrm{d}z$
(using (\ref{eq:SimplifyV})).\vadjust{\goodbreak}

\begin{figure}

\includegraphics{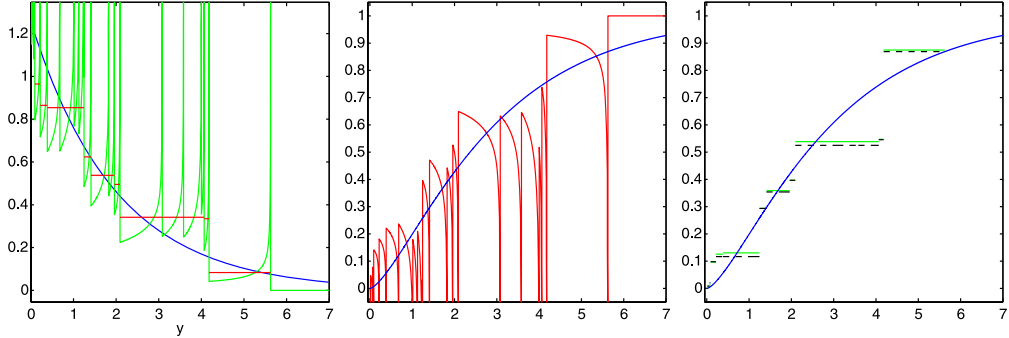}

\caption{Left panel: Plots of $V_n^{\#}$ (green), $\tilde V_n^{\#}$
(red, piece-wise constant) and $V$ (blue, smooth); middle panel:
$F_n^{\#}$ (red), $F$ (blue, smooth); right panel: $\tilde F_n^{\#}$
(green, piece-wise constant), $\check F_n$ (black, piece-wise constant)
and $F$ (blue, smooth) from a sample with $n = 20$ data points.}\label{fig:V_n_F_n}
\end{figure}

Now suppose that we observe an i.i.d. sample $\{(X_{i1},X_{i2})\}
_{i=1}^n$ having the same distribu\-tion as $(X_1,X_2)$. Letting $Y_i =
X_{i1}^2 + X_{i2}^2$, a natural (unbiased) ``naive'' estimator of $V$~is
\[
V_n^{\#}(y) := \int_y^\infty\frac{\mathrm{d}G_n^{\#}(u)}{\sqrt{u - y}} = \frac
{1}{n} \sum_{i=1}^n \frac{\mathbf{1}\{Y_i > y\}}{\sqrt{Y_i - y}},
\]
where $G_n^{\#}$ is the empirical distribution function (EDF) of the
$Y_i$'s. Then $V_n^{\#}(y)$ is an unbiased estimator of $V(y)$ for each
fixed $y$; but $V_n^{\#}$ has infinite discontinuities at the data
points $Y_i$ and is, therefore, not monotonic when viewed as a function
of $y$. See Figure~\ref{fig:V_n_F_n}. We call $V_n^{\#}$ the \textit
{naive estimator}. The naive estimator can be improved by requiring
monotonicity. If $V_n^{\#}$ were square integrable, this could be
accomplished by minimizing the integral of $(W - V_n^{\#})^2$ over all
nonincreasing functions \textit{W}, or equivalently,
\begin{eqnarray}\label{eq:Psi_minimization}
\int_{0}^{\infty} W^2(y) \,\mathrm{d}y - 2 \int_{0}^{\infty} W(y) V_n^{\#}(y) \,\mathrm{d}y.
\end{eqnarray}
The function $V_n^{\#}$ is not square integrable, but it is integrable,
so (\ref{eq:Psi_minimization}) is well defined. Let $\tilde V_n^{\#}$
be the nonincreasing function $W$ that minimizes (\ref
{eq:Psi_minimization}). Existence and uniqueness can be shown along the
lines of Theorem 1.2.1 of Robertson \textit{et~al.} \cite{rwd88}, replacing the sums by
integrals. Groeneboom and Jongbloed \cite{gj95} derived the limit distributions of $V_n^{\#}$
and $\tilde V_n^{\#}$: Let $x_0 > 0$ and
\begin{equation}\label{eq:epsilon_n}
\varepsilon_n := \sqrt{n^{-1} \log n},
\end{equation}
then under appropriate conditions,
\begin{eqnarray}
\label{eq:limitVn-V} {V_n^{\#}(x_0) - V(x_0)\over\varepsilon_n} & \Rightarrow& N (0, g(x_0)
),\\
\label{eq:limitVntilde-V}{\tilde V_n^{\#}(x_0) - V(x_0)\over\varepsilon_n} & \Rightarrow& N \biggl(0,
\frac{1}{2} g(x_0) \biggr),
\end{eqnarray}
where $\Rightarrow$ denotes weak convergence.

We can define two estimators of $F$, $F_n$ and $\tilde F_n^{\#}$, by
replacing $V$ from the right-hand side of (\ref{eq:F}) with $V_n^{\#}$
and $\tilde V_n^{\#}$, respectively.
Observe that $F_n^{\#}$ is not even nondecreasing; $\tilde{F}_n^{\#}$
is nondecreasing, and $\max\{\tilde F_n^{\#}, 0\}$ (as $\tilde F_n^{\#
} \le1$), is a valid distribution function and a more appealing
estimator of $F$ (see Figure~\ref{fig:V_n_F_n}).

Yet another estimator of $F$ can be obtained by isotonizing $F_n^{\#}$
over all nondecreasing functions. Let $\check F_n$ be the nondecreasing
function that is closest to $F_n^{\#}$, in the sense that it minimizes
(\ref{eq:Psi_minimization}) with $V_n^{\#}$ replaced by $F_n^{\#}$. It
is not difficult to see that then $\max\{0,\min(\check F_n, 1)\}$ is a
valid distribution function. Figure~\ref{fig:V_n_F_n} shows the graphs
of the estimators $V_n^{\#}$, $\tilde V_n^{\#}$, $F_n^{\#}$, $\tilde
F_n^{\#}$, and $\check F_n$ obtained from simulated data with $n = 20$.

It will be shown later that for $x_0 > 0$,
\begin{eqnarray}
\label{eq:limitFn-F}\frac{F_n^{\#}(x_0) - F(x_0)}{\varepsilon_n} & \Rightarrow& N \biggl(0, \frac
{4}{\uppi^2} x_0 g(x_0)\biggr ),  \\
\label{eq:limitFntilde-F}\frac{\tilde F_n^{\#}(x_0) - F(x_0)}{\varepsilon_n} & \Rightarrow& N \biggl(0,
\frac{2}{\uppi^2} x_0 g(x_0) \biggr) \quad\mbox{and}  \\
 \label{eq:limitFncheck-F}\frac{\check F_n(x_0) - F(x_0)}{\varepsilon_n} & \Rightarrow& N \biggl(0, \frac
{2}{\uppi^2} x_0 g(x_0) \biggr)
\end{eqnarray}
under modest conditions. As above the isotonized estimators have
exactly half limiting variances of corresponding naive estimators.

Construction of confidence intervals for $F(x_0)$ using these limiting
distributions is still complicated as they require the estimation of
the nuisance parameter $g(x_0)$. Bootstrap intervals avoid this problem
and are generally reliable and accurate in problems with $\sqrt{n}$
convergence rate (see Bickel and Freedman \cite{bf81}, Singh \cite{s81}, Shao and Tu \cite{st95} and its
references). However, conventional bootstrap estimators are
inconsistent for some shape restricted estimators -- dramatically so for
the Grenander estimator. See Kosorok \cite{k08}, Abrevaya and Huang \cite{ah05} and Sen \textit{et~al.} \cite{sbw10}
and its references. So, it is not a priori clear whether bootstrap
methods are consistent in the present context. We show that they are.

In Section~\ref{prelimW}, we prove uniform versions of (\ref
{eq:limitVn-V}), (\ref{eq:limitVntilde-V}), (\ref{eq:limitFn-F}), (\ref
{eq:limitFntilde-F}) and (\ref{eq:limitFncheck-F}). These are used in
Section~\ref{Boots} to establish the consistency of bootstrap methods
in approximating the sampling distribution of the various estimators of
$V$ and $F$, while generating samples from the EDF. Using data on the
globular cluster M62 we illustrate the isotonized estimators of $F$
along with the corresponding bootstrap based point-wise CIs in
Section~\ref{dataana}. Section~A, the Appendix, gives the
details of some of the arguments in the proofs of the main results.
%
\section{Uniform convergence}\label{prelimW}
In this section, we prove central limit theorems for estimates of $V$
and $F$ when we have a triangular array of random variables whose
row-distributions satisfy certain regularity conditions. This
generalization will also help us analyze the asymptotic properties of
the bootstrap estimators (to be introduced in Section~\ref{Boots}).
Note that conditional on the data, bootstrap\vadjust{\goodbreak} samples can be embedded in
a triangular array of random variables, with the $n$th row being
generated from a distribution (built from the first $n$ data points)
that approximates the data-generating mechanism.

Suppose that we have i.i.d. triangular data $\{Y_{n,i}\}_{i=1}^n$
having distribution function $G_n$. We consider a special construction
of $Y_{n,i}$, namely, let $Y_{n,i} = G_n^{-1}(T_{i})$, where
$G_n^{-1}(u) = \inf\{x\dvt G_n(x) \ge u\}$ and $T_{1},T_2,\ldots$ are
i.i.d. $\operatorname{Uniform}(0,1)$ random variables. Let $V_n$ and
$U_n$ be defined as
\begin{eqnarray*}
V_n(y) = \int_y^{\infty}\frac{\mathrm{d}G_n(u)}{\sqrt{u - y}} \quad\mbox{and}\quad
U_n(x) = \int_{0}^{x} V_n (y) \,\mathrm{d}y.
\end{eqnarray*}
Let $\operatorname{LCM}_I$ be the operator that maps a function
$h\dvtx\mathbb{R} \rightarrow\mathbb{R}$ into the least concave majorant
(LCM) of its restriction to the interval $I \subset\mathbb{R}$. Define
$\tilde V_n := \operatorname{LCM}_{[0,\infty)}[U_n]'$ where $'$ denotes
the right derivative. Let $G_n^{\#}$ denote the EDF of
$Y_{n,1},Y_{n,2}, \ldots, Y_{n,n}$,
\[
V_n^{\#}(y) := \int_y^{\infty}\frac{\mathrm{d} G_n^{\#}(u)}{\sqrt{u - y}} = \frac
{1}{n} \sum_{i:Y_{n,i} > y} \frac{1}{\sqrt{Y_{n,i} - y}}.
\]
Then $V_n^{\#}$ is a nonmonotonic, unbiased estimate of $V_n(y)$, as
above, and we call $V_n^{\#}$ the {naive} estimator. The naive
estimator can be improved by imposing the monotonicity constraint as in
(\ref{eq:Psi_minimization}) to obtain $\tilde V_n^{\#}$. Observe that
\[
U_n^{\#}(x) := \frac{2}{n} \sum_{i=1}^n \bigl\{\sqrt{Y_{n,i}} - \sqrt
{(Y_{n,i} - x)_+}\bigr \}
\]
is an unbiased estimate of $U_n(x)$ for all $x \in[0,\infty)$; $U_n^{\#
}$ is a nondecreasing function; $V_n^{\#}$ is the derivative of $U_n^{\#
}$ a.e. Let $\tilde U_n^{\#}$ be the LCM of $U_n^{\#}$. Then $\tilde
V_n^{\#}$ is the right-derivative of $\tilde U_n^{\#}$ (see, e.g.,
Lemma 2 of \cite{gj95}). Let $F_n$ and $F_n^{\#}$ be defined by
replacing $V$ from the right-hand side of (\ref{eq:F}) with $V_n$ and
$V_n^{\#}$, respectively.\vspace*{-2pt}
%
\subsection{CLT for estimates of $V$}
Fix $x_0 \in(0,\infty)$ such that $g(x_0) > 0$. We consider two
estimates of $V(x_0)$, namely $V_n^{\#}(x_0)$ and $\tilde V_n^{\#
}(x_0)$. To find the limit distribution of $V_n^{\#}(x_0)$, we assume
the following conditions on $G_n$:
\begin{eqnarray}
V_n(x_0) & \rightarrow& V(x_0), \label{eq:V_n->V} \\[-1pt]
n \biggl\{G_n \biggl(x_0 + \frac{1}{\varepsilon^2 n \log n}\biggr ) - G_n(x_0)\biggr \} &
\rightarrow& 0 \qquad\mbox{for all } \varepsilon>0, \label{eq:G_n->G} \\[-1pt]
\int_{x_0}^{x_0 + c_n} \frac{\mathrm{d}G_n(y)}{\sqrt{y - x_0}} & = & \mathrm{o}(\varepsilon
_n), \label{eq:IntCond}\\[-1pt]
\frac{1}{\log n} \int_{x_0 + c_n}^\infty\frac{\mathrm{d}G_n(y)}{y - x_0} &
\rightarrow& g(x_0), \label{eq:V_n2->V2}
\end{eqnarray}
where $c_n = 1/(\sqrt{n \log n} + V_n(x_0))^2$ and $\varepsilon_n$ is
defined in (\ref{eq:epsilon_n}).\vadjust{\goodbreak}
\begin{proposition}\label{proposition:limitVn}
If (\ref{eq:V_n->V})--(\ref{eq:V_n2->V2}) hold then $\varepsilon_n^{-1} \{
V_n^{\#}(x_0) - V_n(x_0)\} \Rightarrow N (0, g(x_0) )$.
\end{proposition}

The proof of the proposition is given in the Appendix. Next, we study
the limiting distribution of
\[
\Delta_n := {\tilde V_n^{\#}(x_0) - \hat V_n(x_0)\over\varepsilon_n},
\]
where $\hat V_n(x_0)$ can be $V_n(x_0)$ or $\tilde V_n(x_0)$. Define
the stochastic process
\[
\mathbb{Z}_n(t) = \varepsilon_n^{-2}\{U_n^{\#}(x_0 + \varepsilon_n t) -
U_n^{\#}(x_0) - \hat V_n(x_0) \varepsilon_n t\}
\]
for $t \in I_n:= [-\varepsilon_n^{-1} x_0,\infty)$ and note that $\Delta_n
= \operatorname{LCM}_{I_n}[\mathbb{Z}_n]'(0)$, that is, $\Delta_n$ is
the right-hand slope at $0$ of the LCM of the process $ \mathbb{Z}_n$.
We will study the limiting behavior of the process $\mathbb{Z}_n$ and
use continuous mapping arguments to derive the limiting distribution of
$\Delta_n$. We consider all stochastic processes as random elements in
$C(\mathbb{R})$, the space of continuous functions on $\mathbb{R}$,
equipped it with the Borel $\sigma$-field and the metric of uniform
convergence on compacta. To better understand the limiting behavior of
$\mathbb{Z}_n$, we decompose $\mathbb{Z}_n$ into the sum of
\begin{eqnarray}\label{eq:Simplify_Z}
\mathbb{Z}_{n,1}(t) & = & \varepsilon_n^{-2}\{(U_n^{\#} - U_n)(x_0 +
\varepsilon_n t) - (U_n^{\#} - U_n)(x_0)\} \quad\mbox{and} \nonumber\\
\mathbb{Z}_{n,2}(t) & = & \varepsilon_n^{-2} \{U_n(x_0 + \varepsilon_n t) -
U_n(x_0) - \hat V_n(x_0) \varepsilon_n t\}. \nonumber
\end{eqnarray}
Observe that $\mathbb{Z}_{n,2}$ depends only on $G_n$ and not on the
$Y_{n,j}$. Let
\[
\mathbb{Z}_{1}(t) = t W \quad\mbox{and}\quad \mathbb{Z}(t) = \mathbb{Z}_1(t) +
\tfrac{1}{2} t^2 V'(x_0)
\]
for $t \in\mathbb{R}$, where $W$ is a normal random variable having
mean $0$ and variance $\frac{1}{2} g(x_0)$. We state some conditions on
the behavior of $G_n$, $\hat V_n$ and $U_n$ used to obtain the limiting
distribution of $\Delta_n$.
\begin{itemize}[(b)]
\item[(a)] $D_n:= \| G_n - G \| = \mathrm{O}(\varepsilon_n)$, where $\| \cdot\|$
refers to the uniform norm, that is, $\| G_n - G \| = \sup_{t \in
\mathbb{R}} |G_n(t) - G(t)|$.

\item[(b)] $\mathbb{Z}_{n,2}(t) \rightarrow\frac{1}{2} t^2 V'(x_0)$ as
$n \rightarrow\infty$ uniformly on compacta.

\item[(c)] For each $\varepsilon> 0$,
\[
\bigl| U_n(x_0 + \beta) - U_n(x_0) - \beta\hat V_n (x_0) - \tfrac{1}{2} \beta
^2 V'(x_0) \bigr| \le\varepsilon \beta^2 + \mathrm{o}(\beta^2) + \mathrm{O}(\varepsilon_n^{2})
\]
for large $n$, uniformly in $\beta$ varying over a neighborhood of zero.
\end{itemize}
\begin{them}\label{thm:cons_Z_n}
Under condition \textup{(a)} the distribution of $\mathbb{Z}_{n,1}$ converges to
that of $\mathbb{Z}_1$. Further, if \textup{(b)} holds, then the distribution of
$\mathbb{Z}_n$ converges to that of $\mathbb{Z}$.
\end{them}
\begin{pf} The covariance of $\mathbb{Z}_{n,1}(s)$ and $\mathbb
{Z}_{n,1}(t)$, is needed. To compute it let
\[
\phi(y,\eta) := \sqrt{(y - x_0)_+} - \sqrt{(y - x_0 - \eta)_+}
\]
for $y, \eta\in\mathbb{R}$ and observe the following two properties:
\setcounter{equation}{0}
\renewcommand{\theequation}{P\arabic{equation}}
\begin{eqnarray}
\label{eqp1}\phi(\cdot,\eta) &\le&\sqrt{|\eta|},
\\
\label{eqp2}\int_0^\infty| \phi'(y,\eta)| \,\mathrm{d}y &\le& 2\sqrt{|\eta|}
\end{eqnarray}
of which the second follows from splitting the interval of integration
into $[0,x_0+\eta]$ and $(x_0+\eta,\infty)$. Observe that
\[
\mathbb{Z}_{n,1}(t) = \frac{2}{\varepsilon_n^2} \int\phi(u,\varepsilon_n t)
\,\mathrm{d}(G_n^{\#} - {G}_n)(u)
\]
and
\[
\operatorname{Cov}(\mathbb{Z}_{n,1}(s),\mathbb{Z}_{n,1}(t)) = \frac
{4}{n \varepsilon_n^4} \operatorname{Cov}(\phi(Y_{n,1}, \varepsilon_n s),\phi
(Y_{n,1}, \varepsilon_n t)),
\]
where $Y_{n,1} \sim{G}_n$. We first show that $E[\phi(Y_{n,1}, \varepsilon
_n t)] = \mathrm{O}(\varepsilon_n)$, so that $\operatorname{Cov}(\mathbb
{Z}_{n,1}(s)$, $\mathbb{Z}_{n,1}(t)) = (1/n\varepsilon_n^{4})E[\phi
(Y_{n,1}, \varepsilon_n s)\phi(Y_{n,1}, \varepsilon_n t)] + \mathrm{o}(1)$. For this,
observe that
\setcounter{equation}{4}
\renewcommand{\theequation}{\arabic{section}.\arabic{equation}}
\begin{eqnarray}\label{eq:bound_EZ_n_1}
E[\phi(Y_{n,1}, \varepsilon_n t)] = 2 \int\phi(u,\varepsilon_n t) \,\mathrm{d} ({G}_n
- {G})(u) + \{U(x_0 + t \varepsilon_n) - U(x_0)\}.
\end{eqnarray}
The first term is at most
\begin{eqnarray}\label{eq:bound_Z_n_1}
2\int_0^\infty({G}_n - {G})(u) \phi'(u,\varepsilon_n t) \,\mathrm{d}u & \le& 2 D_n
\int_0^\infty| \phi'(u,\varepsilon_n t)| \,\mathrm{d}u \nonumber\\
& = & 2 \mathrm{O}(\varepsilon_n) 2 \sqrt{|\varepsilon_n t|}=\mathrm{O}(\varepsilon_n^{3/2}),
\nonumber
\end{eqnarray}
and the second term in (\ref{eq:bound_EZ_n_1}) is at most $\mathrm{O}(\varepsilon
_n)$ by using a one term Taylor expansion. Next, suppose that $s \le t$
and write $E[\phi(Y_{n,1},\varepsilon_n s) \phi(Y_{n,1}, \varepsilon_n t)]$ as
\begin{eqnarray}\label{eq:Phi_Expan0}
\int\phi(u, \varepsilon_n s) \phi(u, \varepsilon_n t) \,\mathrm{d} ({G}_n - G)(u) +
\int\phi(u, \varepsilon_n s) \phi(u, \varepsilon_n t) \,\mathrm{d}G(u).
\end{eqnarray}
From Lemma $3$ of Groeneboom and Jongbloed \cite{gj95}, page 1539,
\begin{eqnarray}\label{eq:GJ_result}
\int\phi(u, \varepsilon_n s) \phi(u, \varepsilon_n t) \,\mathrm{d}G(u) = -
\frac{1}{4} g(x_0) s t \varepsilon_n^2 \log\varepsilon_n +
\mathrm{O}(\varepsilon_n^2).
\end{eqnarray}
Using integration by parts, (\ref{eqp1}), and (\ref{eqp2}), the first term in (\ref
{eq:Phi_Expan0}) is at most
\begin{eqnarray}\label{eq:Phi_Expan}
& & \biggl| \int_{0}^{\infty} \{ \phi'(u,\varepsilon_n s) \phi(u,\varepsilon_n t)
+ \phi(u,\varepsilon_n s) \phi'(u,\varepsilon_n t) \} ({G}_n - G)(u) \,\mathrm{d}u \biggr|
\nonumber\\
& & \quad\le2 D_n \sqrt{|\varepsilon_n t|} \bigl\{2 \sqrt{|\varepsilon_n s|}\bigr\} =
\mathrm{O}(\varepsilon_n^2). \nonumber
\end{eqnarray}
So,
\begin{eqnarray}\label{eq:CovZn1}
\operatorname{Cov}(\mathbb{Z}_{n,1}(s),\mathbb{Z}_{n,1}(t)) & = & \frac
{4}{n \varepsilon_n^4} \biggl\{ \frac{1}{4} s t g(x_0) \varepsilon_n^2 \log
\biggl({1\over \varepsilon_n} \biggr) + \mathrm{O}(\varepsilon_n^2) \biggr\} \nonumber
\\[-8pt]
\\[-8pt]
& = & \frac{1}{2} g(x_0) s t \biggl\{1 - \frac{\log\log n}{\log n} \biggr\} + \mathrm{O}
\biggl(\frac{1}{\log n} \biggr).\nonumber
\end{eqnarray}
It follows directly from the Lindeberg--Feller central limit theorem for
triangular arrays that $\mathbb{Z}_{n,1}(1) \Rightarrow N(0,\frac{1}{2}
g(x_0))$; and Chebyshev's inequality implies that $|s \mathbb
{Z}_{n,1}(t) - t \mathbb{Z}_{n,1}(s)| = \mathrm{o}_P(1)$ as $n \rightarrow\infty
$ for all for all fixed $s,t \in\mathbb{R}$. So, the finite
dimensional distributions of $\mathbb{Z}_{n,1}$ converges weakly to the
finite dimensional distributions of $\mathbb{Z}_1$.

For the the convergence in distribution of $\mathbb{Z}_{n,1}$ to
$\mathbb{Z}_{1}$ in $C(\mathbb{R})$, it suffices to show that for each
$M > 0$ and sequence of positive numbers $\{\delta_n\}$ converging to zero,
\[
E\{\sup|\mathbb{Z}_{n,1}(s) - \mathbb{Z}_{n,1}(t)|\dvt |s - t| \le\delta
_n, \max(|s|,|t|) \le M\} \rightarrow0.
\]
See Theorem 2.3 of Kim and Pollard \cite{kp90}. Consider the class of functions
$\mathcal{C}_R = \{\phi(\cdot,\eta) \dvt |\eta| < R\}$ with its natural
envelope $\Phi_R(y):= \sqrt{(y - x_0 + R)_+} - \sqrt{(y - x_0 - R)_+}$.
Observe that $\mathcal{C}_R$ are uniformly manageable for its envelope
$\Phi_R$ and that $\Phi_R \le\sqrt{2R}$. Let $\delta_n$ be a sequence
of positive numbers converging to zero, $h(y;s,t) := \phi(y,t) - \phi
(y,s) = \sqrt{(y - x_0 - s)_+} - \sqrt{(y - x_0 - t)_+}$ for $y,s,t \in
\mathbb{R}$, and $\mathcal{H}_n := \{ h(\cdot;s \varepsilon_n,t\varepsilon
_n)\dvt \max(|s|,|t|) < M, |s - t| \le\delta_n \}$. The class $\mathcal
{H}_n$ has envelope $H_n := 2 \Phi_{M \varepsilon_n}$. Observe that
\begin{eqnarray}
\mathbb{Z}_{n,1}(t) & = & 2 \varepsilon_n^{-2} (G_n^{\#} - G_n) \phi(\cdot
,t \varepsilon_n). \nonumber
\end{eqnarray}
So, it suffices to show that $ \varepsilon_n^{-2} E [ \sup_{h \in\mathcal
{H}_n} |(G_n^{\#} - G_n) h| ]= \mathrm{o}(1)$.
Define $S_n := {G_n^{\#} H_n^2}/{(n \varepsilon_n^4)}$ and $T_n := \sup_{h
\in\mathcal{H}_n} G_n^{\#} h^2$. Then by the maximal
inequality of Section 3.1 in Kim and Pollard \cite{kp90}, there is a (single)
continuous function $J(\cdot)$ for which $J(0) = 0$, $J(1) <\infty$, and
\begin{eqnarray}
{1\over\varepsilon_n^{2}} E\Bigl [ \sup_{h \in\mathcal{H}_n} |G_n^{\#} h -
G_n h| \Bigr]& \le& \frac{1} {\varepsilon_n^{2}\sqrt{n}} E \biggl[ \sqrt{G_n^{\#}
H_n^2} J \biggl(\sup_{\mathcal{H}_n} {G_n^{\#} h^2\over G_n^{\#} H_n^2} \biggr)\biggr ]
\nonumber\\
& = & E \biggl[ \sqrt{S_n} J \biggl( \frac{T_n}{n \varepsilon_n^4 S_n} \biggr) \biggr]. \nonumber
\end{eqnarray}
Let $\eta> 0$. Splitting according to whether $\{S_n \le\eta\}$ or
not, using the fact that $n \varepsilon_n^4 S_n \ge T_n$ and invoking the
Cauchy--Schwarz inequality for the contribution from $\{S_n > \eta\}$,
we may bound the last expected value by
\begin{eqnarray*}
& & E \biggl[ \sqrt{S_n} 1\{S_n \le\eta\} J \biggl(\frac{T_n}{n\varepsilon_n^4 S_n}\biggr )
\biggr] + E \biggl[ \sqrt{S_n} 1\{S_n > \eta\} J \biggl(\frac{T_n}{n \varepsilon_n^4 S_n}\biggr )
\biggr] \\
&&\quad \le \sqrt{\eta} J(1) + \sqrt{ E S_n} \sqrt{E J^2 \biggl( \min\biggl (1,\frac
{T_n}{n \varepsilon_n^4 \eta} \biggr) \biggr)}.
\end{eqnarray*}
Noting that $\Phi_{M \varepsilon_n} = \phi(\cdot,M \varepsilon_n) - \phi(\cdot
, -M \varepsilon_n)$ and using (\ref{eq:GJ_result}) and (\ref{eq:CovZn1})
with $-s = t = M$, we have
\begin{equation}\label{eq:ES_n}
E S_n = \frac{1}{n \varepsilon_n^4} E \Biggl[ \frac{1}{n} \sum_{i=1}^n
H_n^2(Y_{n,i})\Biggr ] = \frac{G_n H_n^2}{n \varepsilon_n^4} = \mathrm{O}(1).
\end{equation}
So, it suffices to show that $T_n = \mathrm{o}_P(n \varepsilon_n^4)$, which implies
$E [ J^2(\min(1,T_n/(n \varepsilon_n^4 \eta)) ] \rightarrow0$ (note that
$J(1) < \infty$). We will establish the stronger result $E T_n = \mathrm{o}(n
\varepsilon_n^4)$. Observe that
\[
E \Bigl[ \sup_{\mathcal{H}_n} G_n^{\#} h^2 \Bigr]  \le E \sup_{\mathcal{H}_n}
G_n h^2 + E \Bigl[ \sup_{\mathcal{H}_n} |G_n^{\#} h^2 - G_n h^2|\Bigr ]
\]
and
\begin{eqnarray}
G_n h^2 &=& G_n [\phi(y,t \varepsilon_n) - \phi(y,s \varepsilon_n)]^2  =
-\tfrac{1}{4} g(x_0) (s - t)^2 \varepsilon_n^2 \log\varepsilon_n + \mathrm{O}(\varepsilon
_n^2) \nonumber\\
& = & \mathrm{O}(\delta_n^2 n \varepsilon_n^4) + \mathrm{O}(\varepsilon_n^2) = \mathrm{o}(n \varepsilon
_n^4) \nonumber
\end{eqnarray}
by (\ref{eq:GJ_result}). The maximal inequality applied to the
uniformly manageable class $\{h^2\dvt h \in\mathcal{H}_n\}$ with envelope
$H_n^2$ bounds the second term by $\tilde J(1)\sqrt{G_n H_n^4/n} \le
8M\varepsilon_n/\sqrt{n} = \mathrm{o}(n \varepsilon_n^4)$, where we have used (\ref
{eq:ES_n}) and the fact that $H_n^2 \le8M\varepsilon_n$. That $\mathbb
{Z}_{n}$ converges in distribution to $\mathbb{Z}$ in $C(\mathbb{R})$
follows directly.
\end{pf}

A rigorous proof of the convergence of $\Delta_n$ involves a little
more than an application of a continuous mapping theorem. The
convergence $\mathbb{Z}_{n} \Rightarrow\mathbb{Z}$ is only in the
sense of the metric of uniform convergence on compacta. A concave
majorant near the origin might be determined by values of the process
long way from the origin; the convergence $\mathbb{Z}_{n} \Rightarrow
\mathbb{Z}$ by itself does not imply the convergence $\operatorname
{LCM}_{I_n} [\mathbb{Z}_{n}] \Rightarrow\operatorname{LCM}_{\mathbb
{R}}[\mathbb{Z}]$. We need to show that $\operatorname{LCM}_{I_n}[
\mathbb{Z}_{n}]$ is determined by values of $\mathbb{Z}_{n}$ for $t$ in
an $\mathrm{O}_P(1)$ neighborhood of the origin. Corollary~\ref
{corollary:localization} shows the convergence of $\Delta_n$, and its
proof is given in the Appendix.
\begin{corollary}\label{corollary:localization}
Under conditions \textup{(a)--(c)}, the distribution of $\Delta_n$ converges
to that of $W \stackrel{d}{=} \operatorname{LCM}_{\mathbb{R}}[\mathbb{Z}]'(0)$.
\end{corollary}
%
\subsection{CLT for estimates of $F$}\label{unifCLTF}
We consider three estimates of $F$, namely $F_n^{\#}$, $\tilde
F_n^{\#}$ and $\check F_n^{\#}$, where $F_n^{\#}$ and $\tilde F_n^{\#}$
are obtained by replacing $V$ from the right-hand side of (\ref{eq:F})
with $V_n^{\#}$ and $\tilde V_n^{\#}$, respectively;
and $\check F_n^{\#}$ is the closest (in the sense of minimizing (\ref
{eq:Psi_minimization}) with $V_n$ replaced with $F_n^{\#}$)
nondecreasing function to $F_n^{\#}$. We start by deriving the limit
distribution of $F_n^{\#}$. Let $\sigma^2 := \operatorname{Var} [\sin
^{-1} \sqrt{1 \wedge(x_0/Y)} ]$ where $Y \sim G$.
\begin{proposition}\label{proposition:limitFn}
If $g(x_0) > 0$ and $\|G_n - G\| \rightarrow0$ as $n \rightarrow\infty
$, then
\begin{eqnarray}\label{eq:limitIntVn-V}
\sqrt{n} \int_{x_0}^\infty\frac{V_n^{\#}(u) - V_n(u)}{2 \sqrt{u}} \,\mathrm{d}u
\Rightarrow N(0, \sigma^2).
\end{eqnarray}
If also (\ref{eq:V_n->V})--(\ref{eq:V_n2->V2}) hold, then
\begin{equation}\label{eq:uniflimitFn-F}
{ F_n^{\#}(x_0) - F_n(x_0)\over\varepsilon_n}  \Rightarrow N \biggl(0, \frac
{4}{\uppi^2} x_0 g(x_0) \biggr).
\end{equation}
\end{proposition}
\begin{pf} For (\ref{eq:limitIntVn-V}), observe that
\begin{eqnarray}
\int_{x_0}^\infty\frac{V_n^{\#}(u)}{2 \sqrt{u}} \,\mathrm{d} u   =   \frac{1}{n}
\sum_{i=1}^n \int_{x_0}^\infty\frac{\mathbf{1}\{Y_{n,i} > u\}}{2\sqrt
{u} \sqrt{Y_{n,i} - u}} \,\mathrm{d}u
 =  \frac{\uppi}{2} - \frac{1}{n} \sum_{i=1}^n \sin^{-1} \sqrt{1 \wedge
\frac{x_0}{Y_{n,i}}} \nonumber
\end{eqnarray}
after some simplification, and (similarly),
\[
\int_{x_0}^\infty\frac{V_n(u)}{2 \sqrt{u}} \,\mathrm{d}u  =  \frac{\uppi}{2} -
\int_0^\infty\sin^{-1} \sqrt{1 \wedge
\frac{x_0}{y}} \,\mathrm{d}G_n(y).
\]
Relation~(\ref{eq:limitIntVn-V}) now follows from the Lindeberg--Feller
CLT. For (\ref{eq:uniflimitFn-F}), first observe that $F_n^{\#}(x_0) -
F_n(x_0)$ may be written as
\[
- \frac{2}{\uppi} \sqrt{x_0} \{V_n^{\#}(x_0) - V_n(x_0)\} - \frac{2}{\uppi
} \int_{x_0}^{\infty}
\frac{V_n^{\#}(u) - V_n(u)}{2 \sqrt{u}} \,\mathrm{d}u.
\]
From Proposition~\ref{proposition:limitVn}, $\varepsilon_n^{-1} \{V_n^{\#
}(x_0) - V_n(x_0)\} \Rightarrow N(0,g(x_0)). $ Relation (\ref
{eq:uniflimitFn-F}) follows directly from this and (\ref{eq:limitIntVn-V}).
\end{pf}

Applying the proposition with $G_n = G$ verifies (\ref{eq:limitFn-F}).
Next, we derive the limiting
distribution of $\tilde F_n^{\#}$.
\begin{proposition}\label{proposition:limitFntilde}
Suppose that \textup{(a)--(c)} hold with $\hat V_n = \tilde V_n$, then,
\[
{\tilde F_n^{\#}(x_0) - \tilde F_n(x_0)\over\varepsilon_n}  \Rightarrow
 N \biggl(0, \frac{2}{\uppi^2} x_0 g(x_0)\biggr ).
\]
\end{proposition}
\begin{pf} As above $\tilde F_n^{\#}(x_0) - \tilde{F}_n(x_0)$ may be
written as
\[\label{eq:breakFn_tilde-F}
\frac{2}{\uppi} \sqrt{x_0} \{ \tilde V_n(x_0) - \tilde{V}_n^{\#
}(x_0) \} + \frac{2}{\uppi} \int_{x_0}^\infty
\frac{\tilde V_n(u) - \tilde V_n^{\#}(u)}{2 \sqrt{u}} \,\mathrm{d}u.
\]
From Corollary~\ref{corollary:localization}, $\varepsilon_n^{-1} \{ \tilde
V_n(x_0) - \tilde V_n^{\#}(x_0) \} \Rightarrow
N (0, {1\over2}{g(x_0)} )$. Integrating by parts, the integral on the
last display is a most
\begin{eqnarray*}
&&\frac{|\tilde U_n(x_0) - \tilde U_n^{\#}(x_0)|}{2 \sqrt{x_0}} + \frac
{1}{4} \biggl| \int_{x_0}^\infty\frac{\tilde
U_n(u) - \tilde U_n^{\#}(u)}{u^{3/2}} \,\mathrm{d}u \biggr| \\
&&\quad \le \frac{\|\tilde U_n - \tilde U_n^{\#} \|}{2 \sqrt{x_0}} + \|
\tilde U_n - \tilde U_n^{\#} \| \frac{1}{2 \sqrt{x_0}} \\
&&\quad=  \frac{\|\tilde U_n - \tilde U_n^{\#} \|}{\sqrt{x_0}} \le\frac{\|
U_n - U_n^{\#} \|}{\sqrt{x_0}}\\
 &&\quad= \mathrm{O}_P(n^{-1/2}) =
\mathrm{o}_P(\varepsilon_n)
\end{eqnarray*}
by Marshall's lemma and maximal inequality 3.1 of Kim and Pollard \cite{kp90} (to bound
$\|U_n - U_n^{\#}\|$). The proposition follows.
\end{pf}

Now let $H_n(x) := \int_0^x F_n(z) \,\mathrm{d}z $ and $H_n^{\#}(x) := \int_0^x
F_n^{\#}(z) \,\mathrm{d}z $. Note that $F_n^{\#}$ is the derivative of $H_n^{\#}$
a.e. Let $\check H_n^{\#}$ be the greatest convex minorant (GCM) of
$H_n^{\#}$. Then $\check F_n$ is the right-derivative of $\check H_n^{\#
}$. We want to study the limit distribution of
\[
\Lambda_n := {\check F_n(x_0) - \hat F_n(x_0)\over\varepsilon_n},
\]
where $\hat F_n$ can be $F_n$ or $\tilde F_n$. Let
\[
\mathbb{X}_n(t) := \varepsilon_n^{-2}\{H_n^{\#}(x_0 + \varepsilon_n t) -
H_n^{\#}(x_0) - \hat F_n(x_0) \varepsilon_n t\}
\]
for $t \in I_n:= [-\varepsilon_n^{-1} x_0,\infty)$. As before, we
decompose $\mathbb{X}_n$ into $\mathbb{X}_{n,1}$ and
$\mathbb{X}_{n,2}$ where
\begin{eqnarray*}\label{eq:Simplify_X}
\mathbb{X}_{n,1}(t) & := & \varepsilon_n^{-2}\{(H_n^{\#} - H_n)(x_0 +
\varepsilon_n t) - (H_n^{\#} - H_n)(x_0)\} \quad\mbox{and}\\
\mathbb{X}_{n,2}(t) & := & \varepsilon_n^{-2} \{H_n(x_0 + \varepsilon_n t) -
H_n(x_0) -\hat F_n(x_0) \varepsilon_n t\}.
\end{eqnarray*}
Let $\operatorname{GCM}_I$ be the operator that maps the restriction of
a function $h\dvtx \mathbb{R} \rightarrow\mathbb{R}$ to the interval $I$
into its GCM, and observe that $\Lambda_n = \operatorname
{GCM}_{I_n}[\mathbb{X}_n]'(0)$. Also let
\[
\mathbb{X}_{1}(t) = t W\quad \mbox{and} \quad\mathbb{X}(t) = \mathbb{X}_1(t) +
\tfrac{1}{2} t^2 f(x_0)
\]
for $t \in\mathbb{R}$,
where $W$ is a normal random variable having mean $0$ and variance
$2x_0 g(x_0)/\uppi^2$ and $f$ is the density of $Z = X_1^2 +X_2^2 +
X_3^2$. The following conditions will be used.
\begin{enumerate}[(b$'$)]
\item[(b$'$)] $\mathbb{X}_{n,2}(t) \rightarrow\frac{1}{2} t^2 f(x_0)$ as
$n \rightarrow\infty$ uniformly on compacta.

\item[(c$'$)] For each $\varepsilon> 0$,
\[
\bigl| H_n(x_0 + \beta) - H_n(x_0) -
\beta\hat F_n (x_0) - \tfrac{1}{2} \beta^2 f(x_0) \bigr| \le\varepsilon
\beta^2 + \mathrm{o}(\beta^2) + \mathrm{O}(\varepsilon_n^{2})
\]
for large $n$, uniformly in
$\beta$ varying over a neighborhood of zero.
\end{enumerate}
\begin{them}\label{thm:cons_X_n}
Under condition \textup{(a)}, the distribution of $\mathbb{X}_{n,1}$ converges
to that of $\mathbb{X}_1$. Further, if \textup{(b$'$)} holds, then the
distribution of $\mathbb{X}_n$ converges to that of $\mathbb{X}$.
\end{them}
\begin{pf} Using the definitions of $F_n^{\#}$, $H_n$ and $H_n^{\#}$,
we may write $H_n^{\#}(x) - H_n(x)$ as
\[
-\frac{2}{\uppi} \biggl[ \int_0^x \sqrt{z} \{V_n^{\#}(z) - V_n(z)\} \,\mathrm{d}z +
\int_0^x \int_{z}^{\infty} \frac{(V_n^{\#} - V_n)(u)}{2 \sqrt{u}} \,\mathrm{d}u \,\mathrm{d}z
\biggr]%
\]
and
\begin{eqnarray}\label{eq:boundVn-V}
\biggl| \int_{z}^{\infty} {(V_n^{\#} - V_n)(u)\over2 \sqrt{u}}\,\mathrm{d}u \biggr| & \le&
\biggl|\int_{z}^{\infty} \int_{z}^{y}
\frac{\mathrm{d}u}{2 \sqrt{y - u} \sqrt{u}} \,\mathrm{d}(G_n^{\#} - G_n)(y) \biggr| \nonumber\\
& = & \Biggl| \int_{z}^{\infty} \Biggl\{ \frac{\uppi}{2} -\sin^{-1}\sqrt{\frac
{z}{y}} \Biggr\} \,\mathrm{d}(G_n^{\#} - G_n)(y) \Biggr| \nonumber\\
& = & \Biggl| \int_{z}^{\infty} \frac{\mathrm{d}}{\mathrm{d}y}\Biggl [ \sin^{-1}\sqrt{\frac{z}{y}} \Biggr]
(G_n^{\#} - G)(y) \,\mathrm{d}y \Biggr| \nonumber\\
& = & \frac{\uppi}{2} \|G_n^{\#} - G_n\| = \mathrm{o}(\varepsilon_n) \qquad\mbox{a.s.},
\nonumber
\end{eqnarray}
using the Law of Iterated Logarithms for $\|G_n^{\#} - G_n\| =
\mathrm{o}(\varepsilon_n)$ a.s. Fix a compact set $K = [-M,M]$. Then
\begin{eqnarray}\label{eq:VnH-Vn}
\mathbb{X}_{n,1}(t) & =& -{2\over\uppi}\int_{x_0}^{x_0+\varepsilon_n t} \frac
{\sqrt{z}[V_n^{\#}(z)-V_n(z)]}{\varepsilon_n^2} \,\mathrm{d}z + \mathrm{o}(1) \nonumber\\
& = & -{2\over\uppi}\sqrt{x_0}\int_{x_0}^{x_0+\varepsilon_n t} \frac{[V_n^{\#
}(z)-V_n(z)]}{\varepsilon_n^2} \,\mathrm{d}z + \mathrm{o}_P(1) \\
& = & -{2\over\uppi}\sqrt{x_0}\mathbb{Z}_{n,1}(t) + \mathrm{o}_P(1) \nonumber
\end{eqnarray}
uniformly on $K$. Note that (\ref{eq:VnH-Vn}) follows as $ | \int
_{x_0}^{x_0+\varepsilon_n t} (\sqrt{z} - \sqrt{x_0})[V_n^{\#
}(z)-V_n(z)]/{\varepsilon_n^2} \,\mathrm{d}z |$ can be bounded, using integration by
parts, by
\begin{eqnarray}\label{eq:IntbyParts}
\bigl|\sqrt{x_0 + \varepsilon_n t} - \sqrt{x_0}\bigr| \max_{|s| \le M} |\mathbb
{Z}_{n,1}(s)| + \max_{|s| \le M} |\mathbb{Z}_{n,1}(s)| \int
_{x_0}^{x_0+\varepsilon_n t} \frac{\mathrm{d}z}{2 \sqrt{z}} = \mathrm{o}_P(1),
\end{eqnarray}
as $ \max_{|s| \le M} |\mathbb{Z}_{n,1}(s)| = \mathrm{O}_P(1)$. The theorem now follows.
\end{pf}
\begin{corollary}\label{corollary:local}
Under conditions \textup{(a)}, \textup{(b$'$)} and \textup{(c$'$)}, the distribution of $\Lambda
_n$ converges to that of $W \stackrel{d}{=}
\operatorname{GCM}_{\mathbb{R}}[\mathbb{X}]'(0)$.
\end{corollary}

The proof is very similar to that Corollary~\ref
{corollary:localization} with the LCMs changed to GCMs. The
modifications are outlined in the Appendix.
%
\section{Consistency of the bootstrap}\label{Boots}
We begin with a brief discussion on the bootstrap. Suppose we have
i.i.d. random variables (vectors) $T_1,T_2,\ldots,T_n$ having an
unknown distribution function $\lambda$ defined on a probability space
$(\Omega,\mathcal{A},P)$ and we seek to estimate the sampling
distribution of the random variable $R_n(\mathbf{T}_n,\lambda)$, based
on the observed data $\mathbf{T}_n =(T_1,T_2,\ldots,T_n)$. Let $\mu_n$
be the distribution function of $R_n(\mathbf{T}_n,\lambda)$. The
bootstrap methodology can be broken into three simple steps:
\begin{enumerate}[]
\item[Step 1:] Construct an estimate $\lambda_n$ of $\lambda$ based on
the data (for example, the EDF).

\item[Step 2:] With $\lambda_n$ fixed, draw a random sample of size $n$
from $\lambda_n$, say $\mathbf{T}_n^* = (T_1^*,T_2^*,\ldots,T_{n}^*)$
(identically distributed and conditionally independent given $\mathbf
{T}_n$). This is called the \textit{bootstrap sample}.

\item[Step 3:] Approximate the sampling distribution of $R_n(\mathbf
{T}_n,\lambda)$ by the sampling distribution of $R^*_n = R_n(\mathbf
{T}_n^*,\lambda_n)$. The sampling distribution of $R^*_n$, the
bootstrap distribution, can be simulated on the computer by drawing a
large number of bootstrap samples and computing $R^*_n$ for each sample.
\end{enumerate}

Thus the bootstrap estimator of the sampling distribution function of
$R_n(\mathbf{T}_n,\lambda)$ is given by $\mu_{n}^*(x) = P^*\{R^*_n \le
x\}$ where $P^*\{\cdot\}$ is the conditional probability given the
data $\mathbf{T}_n$. Let $L$ denote the Levy metric or any other metric
metrizing weak convergence of distribution functions. We say that $\mu
_{n}^*$ is \textit{(weakly) consistent} if $L(\mu_n, \mu_n^*)\stackrel
{P}{\rightarrow} 0$. Similarly, $\mu_{n}^*$ is \textit{strongly
consistent} if $L(\mu_n,\mu_n^*)\stackrel{}{\rightarrow}
0$ a.s. If $\mu_{n}$ has a weak limit $\mu$, for the bootstrap
procedure to be consistent, $\mu_{n}^*$ must converge weakly to $\mu$,
in probability. In addition, if $\mu$ is continuous, we must have
\[
\sup_{x \in\mathbb{R}} |\mu_{n}^*(x) - \mu(x)|
\stackrel{P}{\rightarrow} 0 \qquad\mbox{as } n \rightarrow\infty.\vspace*{0.5pt}
\]
%
\subsection{\texorpdfstring{Bootstrapping $\tilde V_n$}{Bootstrapping V n}}
Given data $Y_1,Y_2,\ldots, Y_n \sim G$ let $G_n^{\#}$ denote its EDF.
Suppose that we draw conditionally independent and identically
distributed random variables $Y_{n,1}^*,Y_{n,2}^*,\ldots,Y_{n,n}^*$
having distribution function $G_n^{\#}$; and let $G_n^{*}$ be the EDF
of the bootstrap sample. Letting
\begin{eqnarray*}\label{eq:Define_Psi_U}
V_n^{*}(y) & := & \frac{1}{n} \sum_{i:Y_{n,i}^* > y} \frac{1}{\sqrt
{Y_{n,i}^* - y}} = \int\frac{\mathbf{1}_{[y,\infty)}(u)}{\sqrt{u - y}}
\,\mathrm{d}G_n^{*}(u) \quad\mbox{and} \\
U_n^{*}(x) & := & \frac{2}{n} \sum_{i=1}^n \bigl\{ \sqrt{Y_{n,i}^*} - \sqrt
{(Y_{n,i}^* - x)_+} \bigr\} \\
& = & 2 \int\bigl\{\sqrt{u} - \sqrt{(u - x)_+}\bigr \} \,\mathrm{d}G_n^{*}
(u),
\end{eqnarray*}
the isotonic estimate of $V$ based on the bootstrap sample is $\tilde
V_n^* = \operatorname{LCM}_{[0,\infty)}[U_n^*]'$. The bootstrap
estimator of the distribution function of $\Delta_n = \varepsilon_n^{-1} \{
\tilde V_n(x_0) - V(x_0)\}$ is then the conditional distribution
function of
$\Delta_n^* := \varepsilon_n^{-1} \{\tilde V_n^*(x_0) - \tilde V_n(x_0)\}$
given the sample $Y_1,\ldots,Y_n$. To find its limit let
\[
\mathbb{Z}_n^*(t) = \varepsilon_n^{-2}\{U_n^{*}(x_0 + \varepsilon_n t) -
U_n^{*}(x_0) - \tilde V_n(x_0) \varepsilon_n t\}
\]
for $t \in I_n:= [-\varepsilon_n^{-1} x_0,\infty)$ and decompose $\mathbb
{Z}_n^*$ into $\mathbb{Z}_{n,1}^*$ and $\mathbb{Z}_{n,2}^*$ where
\begin{eqnarray}\label{eq:Simplify_Z}
\mathbb{Z}_{n,1}^*(t) & = & \varepsilon_n^{-2}\{(U_n^{*} - U_n^{\#})(x_0
+ \varepsilon_n t) - (U_n^{*} - U_n^{\#})(x_0)\}, \nonumber\\
\mathbb{Z}_{n,2}^*(t) & = & \varepsilon_n^{-2} \{U_n^{\#}(x_0 + \varepsilon
_n t) - U_n^{\#}(x_0) - \tilde V_n^{\#}(x_0) \varepsilon_n t\}. \nonumber
\end{eqnarray}
Recall that $\mathbb{Z}_{1}(t) = t W$ and $\mathbb{Z}(t) =
\mathbb{Z}_1(t) + \frac{1}{2} t^2 V'(x_0)$ are two processes defined
for $t \in\mathbb{R}$, where $W$ is a normal random variable having
mean $0$ and variance $\frac{1}{2} g(x_0)$. Let $\mathbf{Y} =
(Y_1,Y_2,\ldots)$. The following theorem shows that bootstrapping
from the EDF $G_n^{\#}$ is weakly consistent.
\begin{them}\label{thm:Cons_Boots_edf}
Suppose that $V$ is continuously differentiable around $x_0$, and
$g(x_0) \ne0$. Then:
\begin{enumerate}[(iii)]
\item[(i)] The conditional distribution of the process
$\mathbb{Z}_{n,1}^*$, given $\mathbf{Y}$, converges to that of
$\mathbb{Z}_1$ a.s.

\item[(ii)] Unconditionally, $\mathbb{Z}_{n,2}^*(t)$ converges in
probability to $\frac{1}{2} t^2 V'(x_0)$, uniformly on compacta.

\item[(iii)] The conditional distribution of the process
$\mathbb{Z}_{n}^*$, given $\mathbf{Y}$, converges to that of
$\mathbb{Z}$, in probability.

\item[(iv)] The bootstrap procedure is weakly consistent, that is, the
conditional distribution of $\Delta_n^*$, given $\mathbf{Y}$,
converges to
that of $W$, in probability.
\end{enumerate}
\end{them}
\begin{pf} Assertion (i) follows directly from Theorem~\ref
{thm:cons_Z_n}, applied with $G_n = G_n^{\#}$, $G_n^{\#}
= G_n^{*}$ and $P\{\cdot\} = P^*\{\cdot\} = P\{\cdot|\mathbf{Y}\}$,
since condition (a) required for Theorem~\ref{thm:cons_Z_n} holds a.s.
For (ii) and (iii), let
\[
\mathbb{Z}_{n}^0(t) = \varepsilon_n^{-2}\{U_n^{\#}(x_0 + t \varepsilon_n) -
U_n^{\#}(x_0) - \varepsilon_n t V(x_0)\}
\]
for $t \in I_n$. By Theorem~\ref{thm:cons_Z_n}, applied with $G_n = G$,
$V_n =
V$ and $U_n = U$ for all $n$, $\mathbb{Z}_{n}^0$ converges in
distribution to $\mathbb{Z}$. To
prove (ii) observe that
\[
\mathbb{Z}_{n,2}^*(t) = \mathbb{Z}_{n}^0(t) - t \cdot\operatorname
{LCM}_{I_n} [\mathbb{Z}_{n}^0]'(0).
\]
Unconditionally, using the continuous mapping theorem along with a
localization argument as in Corollary~\ref{corollary:localization}, we
obtain $\mathbb{Z}_{n,2}^*(t) \Rightarrow\mathbb{Z}(t) - t \cdot
\operatorname{LCM}_{\mathbb{R}}[\mathbb{Z}]'(0) \stackrel{}{=} \frac
{1}{2} t^2 V'(x_0)$. As the limiting process is a constant, $\mathbb
{Z}_{n,2}^*(t) \stackrel{P}{\rightarrow} \frac{1}{2} t^2 V'(x_0)$. Let
$\{n_k\}$ be a subsequence of $\mathbb{N}$. We will show that there
exists a further subsequence such that conditional on
$\mathbf{Y}$, $\mathbb{Z}_n \Rightarrow\mathbb{Z}$ a.s. along the
subsequence. Now, given $\{n_k\}$, there exists a further subsequence
$\{n_{k_l}\}$ such that $\mathbb{Z}_{n_{k_l},2}^*(t) \stackrel
{}{\rightarrow} \frac{1}{2} t^2 V'(x_0)$ uniformly on
compacta a.s. Thus, the conditional distribution of $\mathbb
{Z}_{n_{k_l}}^*$ given $\mathbf{Y}$, converges to that of
$\mathbb{Z}$, for a.e. $\mathbf{Y}$. This completes the proof of (iii).

For (iv), we use Corollary~\ref{corollary:localization}. Although
conditions (a) and (b) hold in probability, condition~(c) holds with
$\hat V_n = \tilde V_n$ and the $\mathrm{O}(\varepsilon_n^{2})$ term replaced by
$\mathrm{O}_P(\varepsilon_n^{2})$. Thus we cannot appeal directly to Corollary~\ref
{corollary:localization}. Let $\xi> 0$ and $\eta> 0$ be given. We
will show that there exists $N \in\mathbb{N}$ such that for all $n \ge
N$, $P\{L(K, K_n^*) > \xi\} < \eta$, where $L$ is the Levy
metric\vadjust{\goodbreak}
(Gnedenko and Kolmogorov \cite{gk68}, page 33), $K$ is the distribution function of $W \sim
N(0,\frac{1}{2}g(x_0))$ and $K_n^*$ is the distribution function of
$\Delta_n^*$, conditional on the
data. For $\varepsilon> 0$, sufficiently small, let
\[
A_n := \bigl\{ \bigl| U_n^{\#} (x_0 + \beta) - U_n^{\#}(x_0) - \beta\tilde
V_n^{\#}(x_0) - \tfrac{1}{2} \beta^2 V'(x_0)\bigr | < C \varepsilon_n^2 +
\varepsilon\beta^2\bigr \},
\]
where $C > 0$ is chosen such that $P\{A_n^c\} < \frac{\eta}{2}$. This
can be done since (c) holds with $\mathrm{O}(\varepsilon_n^{2})$ term replaced by
$\mathrm{O}_P(\varepsilon_n^{2})$. Further, let $P_n^{0}(E)(\omega) = P^*(E)(\omega
)$, if $\omega\in A_n$ and $P_n^{0}(E)(\omega) = P(E)$, if $\omega
\notin A_n$; and let $K_n^0$ be the distribution function of $\Delta
_n^*$ under the probability measure $P_n^0$. Observe that $K_n^{0} =
K_n^{*}$ on $A_n$ and that $L(K, K_n^0) \stackrel{P_n^{0}}{\rightarrow}
0$ by Corollary~\ref{corollary:localization} can be applied. Therefore,
for all sufficiently large $n$,
\begin{eqnarray*}
P\{L(K, K_n^*) > \xi\} &\le& P \biggl\{ L(K, K_n^0) > \frac{\xi}{2} \biggr\} + P
\biggl\{ L(K_n^0, K_n^*) > \frac{\xi}{2} \biggr\} \\
& \le& \frac{\eta}{2} + P\biggl \{ L(K_n^0, K_n^*) > \frac{\xi}{2} , A_n^c
\biggr\} \le \eta.
\end{eqnarray*}
This completes the proof of (iv).
\end{pf}
\begin{rem} Let $J_n(t) = \int_t^\infty\sqrt{z - t} \,\mathrm{d}V_n^{\#}(z)$, for
$t \ge0$, as in (\ref{eq:Simplify_G}). Then ${G}^{\#}_n = 1 +
2J_n(t)/\uppi$ after some simplification. So, using (\ref{eq:Simplify_G})
to generate the bootstrap sample would lead back to ${G}^{\#}_n$.
\end{rem}
%
\subsection{\texorpdfstring{Bootstrapping $F_n$, $\tilde F_n$ and $\check F_n$}{Bootstrapping F n, F n and F n}}
Bootstrap versions of the three estimators of $F$ under study,
$F_n^{*}$, $\tilde F_n^{*}$ and $\check F_n^{*}$ say, are defined as in
Section~\ref{unifCLTF}; for example, $F_n^*(x) = 1 +
(2/\uppi) \int_x^\infty\sqrt{z} \,\mathrm{d}V_n^{*}(z)$. We approximate the
sampling distribution of $\varepsilon_n^{-1}\{F_n^{\#}(x_0) - F(x_0)\}$
by the bootstrap distribution of $\varepsilon_n^{-1}\{F_n^{*}(x_0) -
F_n^{\#}(x_0)\}$. The bootstrap samples are generated from
$G_n^{\#}$, the EDF of the $Y_i$'s. By appealing to Proposition~\ref
{proposition:limitFn} with $G_n = G_n^{\#}$, it is easy to see that the
bootstrap method is weakly consistent as (\ref{eq:V_n->V})--(\ref
{eq:V_n2->V2}) hold in probability.

The sampling distribution of $\varepsilon_n^{-1}\{\tilde F_n^{\#}(x_0) -
F(x_0)\}$ is approximated by that of $\varepsilon_n^{-1}\{\tilde
F_n^{*}(x_0) - \tilde F_n(x_0)\}$. Using Proposition~\ref
{proposition:limitFntilde}, we can establish the consistency of the
method. Note that the proof of Theorem~\ref{thm:Cons_Boots_edf} shows
how conditions (a)--(c) are satisfied with $G_n = G_n^{\#}, \hat V_n =
\tilde V_n$ required to apply Proposition~\ref{proposition:limitFntilde}.

Recall that $\check F_n^*$ is the nondecreasing function closest to
$F_n^{*}$. Let $H_n^{\#}(x) := \int_0^x F_n^{\#}(z) \,\mathrm{d}z $ and
$H_n^{*}(x) := \int_0^x F_n^{*}(z) \,\mathrm{d}z $. Next, we show that
approximating the distribution of $\Lambda_n = \varepsilon_n^{-1} \{\check
F_n(x_0) - F(x_0)\}$ by the bootstrap distribution of $\Lambda_n^* :=
\varepsilon_n^{-1} \{\check F_n^*(x_0) - \check F_n(x_0)\}$ is consistent.
To find the limit of the conditional distribution of $\Lambda_n^*$, let
\[
\mathbb{X}_n^*(t) = \varepsilon_n^{-2}\{H_n^{*}(x_0 + \varepsilon_n t) -
H_n^{*}(x_0) - \check F_n(x_0) \varepsilon_n t\}
\]
for $t \in I_n:= [-\varepsilon_n^{-1} x_0,\infty)$ and decompose it into
$\mathbb{X}_{n,1}^*$ and $\mathbb{X}_{n,2}^*$, where
\begin{eqnarray*}\label{eq:Simplify_X}
\mathbb{X}_{n,1}^*(t) & = & \varepsilon_n^{-2}\{(H_n^{*} - H_n^{\#})(x_0
+ \varepsilon_n t) - (H_n^{*} - H_n^{\#})(x_0)\}, \\
\mathbb{X}_{n,2}^*(t) & = & \varepsilon_n^{-2} \{H_n^{\#}(x_0 + \varepsilon
_n t) - H_n^{\#}(x_0) - \check F_n(x_0) \varepsilon_n t\}.
\end{eqnarray*}
Recall that $\mathbb{X}_{1}(t) = t W$ and $\mathbb{X}(t) = \mathbb
{X}_1(t) + \frac{1}{2} t^2 f(x_0)$ are two processes defined
for $t \in\mathbb{R}$, where $W$ is a normal random variable having
mean $0$ and variance $\frac{2}{\uppi^2} x_0 g(x_0)$.
\begin{them}\label{thm:Cons_Boots_F}
Suppose that $F$ is continuously differentiable around $x_0$, and
$g(x_0) \ne0$. Then:
\begin{enumerate}[(iii)]
\item[(i)] The conditional distribution of the process
$\mathbb{X}_{n,1}^*$, given $\mathbf{Y}$, converges to that of
$\mathbb{X}_1$ a.s.

\item[(ii)] Unconditionally, $\mathbb{X}_{n,2}^*(t)$ converges in
probability to $\frac{1}{2} t^2 f(x_0)$, uniformly on compacta.

\item[(iii)] The conditional distribution of the process
$\mathbb{X}_{n}^*$, given $\mathbf{Y}$, converges to that of
$\mathbb{X}$, in probability.

\item[(iv)] The bootstrap procedure is weakly consistent, that is, the
conditional distribution of $\Lambda_n^*$, given $\mathbf{Y}$,
converges to that of $W$, in probability.
\end{enumerate}
\end{them}
\begin{pf} The proof is very similar to that of Theorem~\ref
{thm:Cons_Boots_edf}. To find the conditional distribution of $\mathbb
{X}_{n,1}^*$ given $\mathbf{Y}$, we appeal to Theorem~\ref
{thm:cons_X_n} with $G_n = G_n^{\#}$, $G_n^{\#} = G_n^{*}$ and $P\{\cdot
\} = P^*\{\cdot\} = P\{\cdot| \mathbf{Y}\}$. Note that condition (a)
required for Theorem~\ref{thm:cons_X_n} holds a.s. We express $\mathbb
{X}_{n,2}^*(t)$ as $\mathbb{X}_{n}^0(t) - t \cdot\operatorname
{GCM}_{I_n} [\mathbb{X}_{n}^0]'(0)$ where
\[
\mathbb{X}_{n}^0(t) = \varepsilon_n^{-2} \{H_n^{\#}(x_0 + t \varepsilon_n) -
H_n^{\#}(x_0) - F(x_0) \varepsilon_n t\}.
\]
Note that unconditionally $\mathbb{X}_{n}^0$ converges in distribution
to $\mathbb{X}$ by an application of Theorem~\ref{thm:cons_X_n} with
$G_n = G$, $\hat F_n = F$ and $H_n = H$ for all $n$.

Unconditionally, using the continuous mapping theorem along with a
localization argument as in Corollary~\ref{corollary:localization}, we
obtain $\mathbb{X}_{n,2}^*(t) \Rightarrow\mathbb{X}(t) - t \cdot
\operatorname{GCM}_{\mathbb{R}}[\mathbb{X}]'(0) \stackrel{}{=} \frac
{1}{2} t^2 f(x_0)$. As the limiting process is a constant, $\mathbb
{X}_{n,2}^*(t) \stackrel{P}{\rightarrow} \frac{1}{2} t^2 f(x_0)$.

An argument using subsequences as in the proof of (iii) of
Theorem~\ref{thm:Cons_Boots_edf} shows that the conditional
distribution of the process $\mathbb{X}_{n}^*$, given $\mathbf{Y}$,
converges to that of $\mathbb{X}$, in probability. The last part of the
theorem follows along similar lines as in the proof of (iv) of
Theorem~\ref{thm:Cons_Boots_edf}.
\end{pf}
%

\begin{figure}

\includegraphics{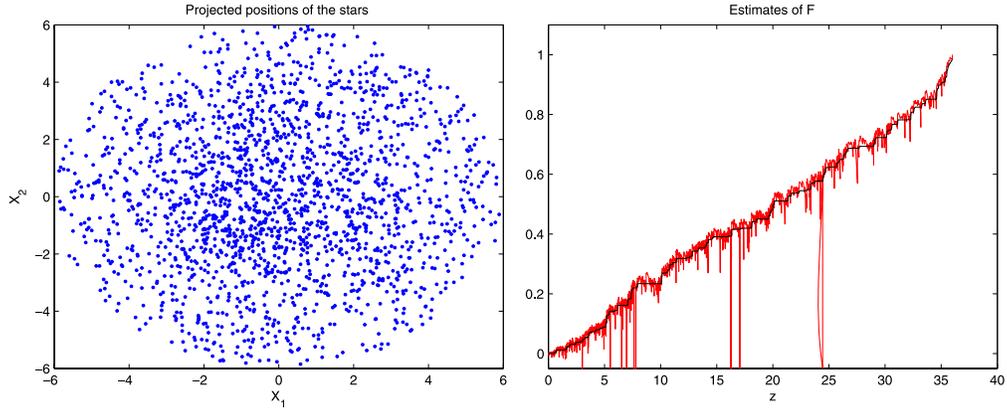}

\caption{Left panel: Scatter plots of the projected positions of the
stars; right panel: $F_n^{\#}$ (red), $\tilde F_n^{\#}$ (black,
piece-wise constant).}\label{fig:Scatterplot}
\end{figure}
%
\begin{figure}[b]

\includegraphics{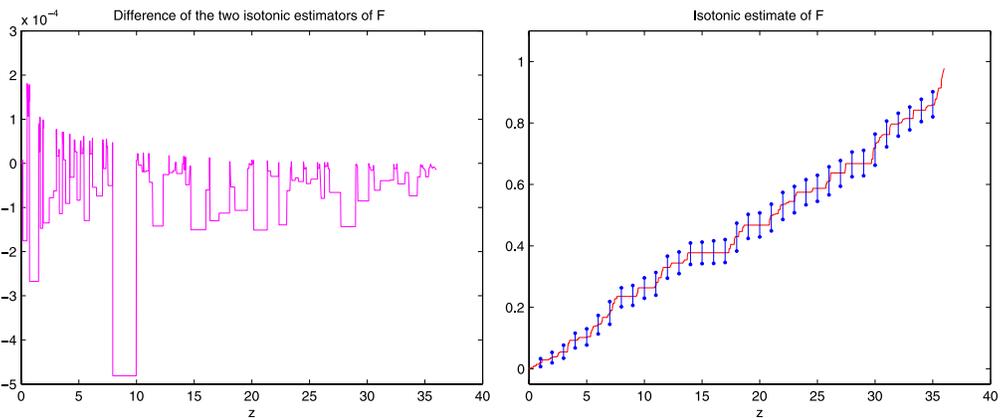}

\caption{Left panel: $\tilde F_n^{\#} - \check{F}_n$; right panel:
Bootstrap based 95\% pointwise CIs around $\tilde F_n^{\#}$.}\label{fig:BootsCIplot}
\end{figure}
%

\section{Data application}\label{dataana}
A globular cluster (GC) is a spherical collection of stars that orbits
a galactic core as a satellite. GCs are very tightly bound by gravity,
which gives them their spherical shapes and relatively high stellar
densities toward their centers. The study of the inner Galactic GCs is
important for several reasons -- to understand the morphology of the
inner Galaxy, to better constraint the characteristics of the Galactic
bulge, etc. Data is available on individual stars in 25 globular
clusters located toward the center of the Milky Way (see, e.g., Alonso \cite
{a10} and Alonso \textit{et~al.} \cite{ams07}). The left panel of Figure \ref
{fig:Scatterplot} shows the projected positions of $n = 2000$ stars in
the inner core of the globular cluster M62 (also known as NGC 6266).
Interest focuses on estimating the distribution function $F$ of the
squared radius. The naive estimator of $F$, $F_n^{\#}$, is shown in the
right panel of Figure~\ref{fig:Scatterplot} along with the isotonized
estimator $\tilde F_n^{\#}$. The two isotonic estimators $\tilde F_n^{\#
}$ and $\check F_n$ are virtually indistinguishable, and the left panel
of Figure~\ref{fig:BootsCIplot} shows the difference between the two
estimators. Note that both the isotonic estimators have the same
pointwise normal limit distribution. The right panel of Figure \ref
{fig:BootsCIplot} shows the point-wise bootstrap based 95\% CIs for
$F$ using the estimator $\tilde F_n^{\#}$. A very similar plot is
obtained using the estimator $\check F_n$.

\section*{Acknowledgements}
The first author's research was supported by the National Science
Foundation, USA.

\begin{supplement}
\stitle{Proofs}
\slink[doi]{10.3150/12-BEJ378SUPP}
\sdatatype{.pdf}
\sfilename{bej378\_supp.pdf}
\sdescription{The Appendix gives the details of some of the arguments
in the proofs of the
main results.}
\end{supplement}


\printhistory

\end{document}